\documentclass{amsart}
\title{The Kolmogorov Reform of Mathematics Education in the USSR}
\author{Alexandre Borovik}
\address{University of Manchester, United Kingdom}
\email{alexandre$\gg$at$\ll$lborovik.net}
\date{11 January 2020}
\thanks{This is a pre-print of the following chapter: Alexandre Borovik, The Kolmogorov Reform of Mathematics Education in the USSR, published in \textbf{Modern mathematics. An international movement?}, edited by Dirk De Bock, 202X, Springer, reproduced with permission of Springer Nature Switzerland AG. The final authenticated version will become available online at: http://dx.doi.org/[to be confirmed].}

\usepackage{hyperref}
\usepackage{amsmath}

\newcommand{\bi}{\begin{itemize}}
\newcommand{\ei}{\end{itemize}}
\newcommand{\bd}{\begin{description}}
\newcommand{\ed}{\end{description}}
\newcommand{\bq}{\begin{quote}
\newcommand{\eq}{\end{quote}}}

\begin{document}


\begin{abstract}
  In the Soviet Union a reform movement in mathematics education was triggered by Andrey Kolmogorov in the 1970s, and followed by a counter-reform. This movement was rooted in the very different socioeconomic conditions of that time and place, and followed a strategy with very significant contrasts to similar programs in the USA, England, or France. This provides an interesting case study which may illuminate the way such movements arise and succeed or fail, and, at the social level, certain fundamental commonalities of constraints as well as significant differences according to local conditions. We shall show that the principal reasons of the failure of the Kolmogorov reform were political: (1) The reform ignored the reality of the socio-economic conditions of the country; (2) The human factor was ignored, and very little attention was given to professional development and retraining of, and methodological help to, the whole army of teachers; (3) An attempt to transfer mathematical content and methods from the highly successful advanced extension stream for mathematically strong and highly engaged children to mainstream education was an especially grievous error.
\end{abstract}
\maketitle


\section*{Introduction}

It is now widely accepted (I have never seen or heard claims to the contrary), that the reform of school mathematics education in the Soviet Union in the 1970s (the reform), initiated and led by the great mathematician Andrey Kolmogorov (1903–1988) was a fiasco---this view has been shared both by supporters of Kolmogorov and by his direct enemies.

What is still disputed is the attribution of guilt. Who was responsible for the damage to the system of mathematics education which, as many of my colleagues in Russia feel, still has not been repaired? The reformers? Bureaucrats from the Ministry of Education? School teachers? I do not wish to take sides in this dispute, but political infighting continues, and one side of this polemic, squarely blaming the reformers, can be seen, for example, in Kostenko (2013, 2014); his works also demonstrate the biggest black spot in the study of this episode in history: The lack of access to archive documents from the highest echelons of power of that time: The Central Committee of the Communist Party, the Ministry of Education of the Soviet Union, and the Ministry of Education of the Russian Federative Republic of the Soviet Union. In the monolithic administrative structure of a totalitarian state, publicly available materials show only the tip of the iceberg of the actual decision making process in the late 1960s which produced the reform launched in 1970.

A rebellion against the reform in the years 1978–1980 appears to be better documented and better understood (Kolyagin \& Savvina, 2012). This rebellion started in the open with the meeting of the Mathematics Division of the Academy of Sciences on December 5, 1978. The symbolic moment marking the beginning of the end of the reform was the publication of a paper by Lev Pontryagin, a famous mathematician, in \emph{Kommunist}, the all-important political journal of the Central Committee of the Communist Party (Pontryagin, 1980). The publication of a paper in \emph{Kommunist} was the seal of approval of the author’s position given from the highest levels of the Party’s hierarchy.

Pontryagin (1980) begins his paper by quoting the definition of vector given in a reformist school textbook by Klopsky, Skopets, and Yagodovsky (1980):
\bq
The vector (parallel translation) defined by a pair $(A, B)$ of distinct points is the transformation of the space which sends each point $M$ to the point $M_1$ such that the ray $MM_1$ is co-directed with the ray $AB$ and the distance $|MM_1|$ equals $|AB|$. (Klopsky et al., 1980, p. 42, quoted in Pontryagin, 1980, p. 99)
\eq
Pontryagin continues as follows:
\bq
This tangle of words is difficult to sort out, but, importantly---it is useless, since it cannot be applied neither in physics, nor in mechanics, nor in other sciences.

What is this? Mockery? Or unintentional absurdity? No, the replacement in these textbooks of many relatively simple, visual formulations by cumbersome, deliberately complicated ones, it turns out, is motivated by the desire... to improve (!) the teaching of mathematics. If the example I gave was just an annoying exception, then the mistake would be easy to eliminate. But, in my opinion, unfortunately, the whole system of school mathematics education has come to a similar state… (p. 99)
\eq
In the political tradition of the Soviet Union, the use of this kind of language was equivalent to branding with hot iron.
The failure of the reform was surprising on a number of counts. First, it was run by mathematicians and professionals in mathematics education of the highest class. Second, The New Math movement in the West was looked at but not imitated. The Kolmogorov reform was not, strictly speaking, the “New Math”---the nature of the reform and the socio-political environment were very different from that of, say, the United States of America (see, e.g., Phillips, 2015). Third, the scope of the reform was relatively modest. Set-theoretic concepts were used but did not play a leading role in the exposition---but nonetheless, the set-theoretic language was used sometimes in a rather annoying way---for example, for introducing central concepts in the course of geometry, such as vector, without further use. No abstract algebra was introduced; the treatment of algebra was mostly untouched but was expanded by inclusion of some elementary calculus. In geometry, the principal change was the introduction of vectors and the systematic use of geometric transformations; again, nothing special, at first glance---this could easily have been done without any mention of set-theoretic concepts.

There was also another aspect. As Sharygin (2002) formulated it:
\bq
It is interesting that the Soviet system of work with mathematically gifted children, created by disinterested enthusiasts and brought, oddly enough, to the level of “know-how,” turned out to be almost the only market product of the Russian education system (not counting, of course, its final result---scientists). (no pagination)
\eq
The key players of the reform: Vladimir Boltyansky (1925–2019), Aleksey Markushevich (1908–1979), Naum Vilenkin (1920–1991), and Isaak Yaglom (1921–1988) happened to be major contributors to it, and moreover, Kolmogorov---one of the creators of that very “Soviet system of work with mathematically gifted children,” as Sharygin highlighted. Surprisingly, there is no umbrella term for the plethora of activities involved. I suggest to use the words “Olympiad Stream” because of its historically oldest and central component: mathematical “Olympiads,” competitions for schoolchildren. Correspondingly, people who passed through the Olympiad Stream and were shaped by it are known as olympiadniks.

I myself was a stereotypical olympiadnik, and, observing the reform unfolding (without affecting me) and collapsing, was puzzled how people whose names were known to me and my friends, whose books we read and respected, managed to botch the reform so spectacularly. In this chapter, I am trying to analyze the intriguing episode in history of mathematics education from this peculiar point of view, that of a mathematician schooled in Russia during the reform period, which, I hope, deserves some attention---simply because it gives a perspective different from the ones used before.

A short contribution like this one allows to mention only a few key actors in the events, and only a few textbooks. I skip entirely the “counter-reform” movement which started in 1980, which also was very important and interesting, with many prominent mathematicians and experienced educationalists involved. Also, I do not make any comparisons with “New Math” reforms in the West---the socio-political environment of the Kolmogorov reform was so different, that comparisons simply make no sense.

A few words on existing literature are due. Karp and Vogeli (2010) provide a good general survey of Russian mathematics education. Neretin (2019) perhaps is the best available analysis of Kolmogorov’s reform. It is 80 pages long, full of detail, and is written from a measured, rational, and non-political position. It also contains a massive bibliography and useful biographic notes. Of other sources, I would recommend Verner (2012), written by an active participant of the “counter-reform.” Kolyagin and Savvina (2012) provide a number of important documents and a well-written introduction with the summary of events. Kolyagin was an active reformer, and, in particular, participated in writing books for teachers (the area of the reform which was neglected by the reformers): Gusev, Kolyagin, and Lukankin (1976), and Kolyagin (1977). Perhaps a close contact with teachers inevitable in writing these books led to Kolyagin’s conversion: He became one of the early critics of the reform. In particular, he was invited to speak at the Meeting of the Mathematics Division of the Academy of Sciences in December 1978 which could now be seen as the start of the “rebellion” of the Academy of Sciences and the Ministry of Education against the reform (Kolyagin \& Savvina, 2012). Kostenko (2013, 2014) is excessively politicized to my taste but is a useful source of bibliographic references. Abramov (2016) is an example of unreserved praise to Kolmogorov.

\section*{Prehistory: The Start of the Olympiad Movement in the 1930s}

To understand what happened in the reform we have to take a look at the early origins of the Olympiad Stream, more specifically, the olympiad movement in the 1930s. A good source for the early history of the movement is Boltyansky and Yaglom (1965).

The first mathematical olympiads were organized in 1934 in Leningrad by Boris Delaunay, a well-known mathematician (Chistyakov, 1935; Fomin, 2020), in 1935 in Moscow by Pavel Alexandrov and other Moscow mathematicians (Bonchkovsky, 1936), and in 1936 in Kazan by the well-known algebraist Nikolai Chebotarev (Chebotarev, 1937). The problems of the Kazan Olympiad were quite challenging (Grigoriev, 1937), which suggested an academically selective approach.

It is interesting that Delaunay’s motives for organization of the olympiad were political: He arranged for young people who were successful in the competition to be admitted to the Mathematics Department of the Leningrad University without formal entrance examinations. This allowed to save some aspiring young mathematicians of “wrong” social backgrounds, that is, children of people deprived of citizens’ rights: army officers and civil servants of the previous tsarist regime, clergy, nobility, capitalists, etc. (Zalgaller, 2021). In the 1920s and early 1930s, it was forbidden to admit children of “lishentsy” [“deprived”] to universities---the olympiad was a loophole.

Mathematicians invested remarkable energy and effort in this project. Why? Because olympiads and other outreach activities were giving them some influence on who is coming into mathematics, and gave a chance to protect university mathematics from political appointees---the latter were quite prominent, for example, in biology---Trofim Lysenko being the most notorious case (Strunnikov \& Shamin, 1989).

In the first Moscow mathematics olympiad of 1935, Pavel Alexandrov was the chairman of the organizing committee. Next year, he wrote in the introduction to a little booklet with problems and solutions of this olympiad (Bonchkovsky, 1936):
\bq
The olympiad is the first entry of future mathematicians into the mathematical arena. It should help us to select these future mathematicians from among our youth, it should help us to provide opportunities for their further mathematical development and education. (p. 4)
\eq
Here we see an unashamed elitism in a supposedly egalitarian country---but this is not so surprising. What is really astonishing is the phrase “It should help us to select these future mathematicians from among our youth.” (p. 4, italics added).

As explained in Borovik, Karakozov, and Polikarpov (2021), selection and promotion of the “cadre” was the ultimate monopoly of the ruling Communist Party. The speech by Andrei Bubnov, People’s Commissary of Education in 1929–1937 at the 17th Party Congress (VKP(b), 1934) is quite illuminating. Besides what would now be described as “widening participation,” ensuring the steady progress of working class children through the school system, he also emphasizes a different task which could be formulated as educating the new generation of loyal to the regime high level specialists for the military, industry, science, medicine. Some examples given by Bubnov in his speech:
\bq
Look at Kamai---professor at Kazan University, Tatar, former docker, now the author of a number of research works in the area of organic compounds of phosphorous and arsenic. (p. 114)
\eq
We should not forget that “organic compounds of phosphorous and arsenic” were an obvious euphemism for “precursors of chemical weapons.” Bubnov’s previous post was that of Head of the Main Political Directorate of the Red Army.

The young mathematicians Pavel Alexandrov, Boris Delaunay, Alexander Gelfond, Alexander Khinchin, Andrey Kolmogorov, Lazar Lyusternik, and Lev\linebreak Schnirelmann, offered, at the right moment, their services to the Party, thus ensuring some degree of their own control over the supply of fresh blood to the top tier of the mathematical profession. In mathematics, “Red Professors,” recruited from the working class party activists (or Young Communist League activists), could not compete with people who had a deeper education because they started their development as mathematicians within the Olympiad Stream.

Mathematicians dared to ask for autonomy in selection and development of their own. It could be seen in other documents of the epoch, for example, in Resolutions of the Second All-Union Congress of Mathematicians which took place in 1934: A special resolution was about olympiads (VSM, 1935, p. 56), and it called for “identification of gifted youths,” stating that “Universities might use olympiads for recruitment of students to mathematical, mechanical and physics departments” and requested funding from the Narkompros (the Ministry of Education) for running olympiads and related activities.

Apparently, these requests were met by the authorities. Starting from 1935, mathematical olympiads and associated activities, first of all, mathematical circles, flourished in Leningrad and Moscow, with crème de la crème of research mathematicians actively involved; what is important, new didactic approaches to exposition of non-trivial mathematics were invented and tested. A good and detailed narrative of these remarkable developments can be found in Boltyansky and Yaglom (1965).

There were two dramatic episodes in the 1930s which also helped mathematicians to gain certain autonomy in running their academic affairs and controlling the intake into the professional mathematical community. One of them was the political, by its nature, struggle around the quality of mathematical textbooks for schools: It is analyzed in Borovik, Karakozov, and Polikarpov (2021); some details in this section are borrowed from that paper. Another one was Luzin’s Afair (Demidov \& Levshin, 1999; Neretin, 2021a). Pavel Alexandrov and other young mathematicians launched a political attack at Nikolai Luzin, a prominent mathematician and the teacher of many of them. I will argue in another paper in preparation, that one of the reasons for the surprisingly vicious attack was Luzin’s tendency to write laudatory reference letters, including letters for political appointees---potential “Red Professors,” which undermined his younger colleagues’ fight for the control of mathematics.

Some may say that establishing this special relationship with the totalitarian regime was entering into a pact with the Devil. I do not want to be judgemental---what mattered, the leading mathematicians established themselves as guardians of the quality of mathematics education and also of mathematical culture; this term is interesting, and, I think, its specific use in Russia is not widely known abroad. It could be traced back at least to 1941, when Otto Schmidt, a mathematician, a member of the Academy of Sciences, and a high ranking government official, formulated the role of mathematicians as specialists who maintain the strategically important mathematical culture of the country:
\bq
Not only us, professionals of science, but the whole country was happy to learn about the solution of a problem set 150 years ago. This problem was solved by academician Vinogradov who proved that every odd number could be written as a sum of three prime numbers.\footnote{This claim was exaggerated: Vinogradov' result was weaker than the one
formulated by Schmidt, although still very impressive (Vavilov, 2021).} Is this needed at the practical level? No. Maybe it is not needed at all? On the contrary, it is much needed, because the culture, the mathematical culture depends on the level of these works in pure mathematics and theoretical physics. You all know that this is not needed for each of us, but we all are interested in the highest possible level of mathematical culture in our country, because it is important to be able to solve any mathematical problem and for that it is important to be able to solve problems such as Goldbach’s problem. The level of mathematical culture in our country is exceptionally high. One may confidently say that in that respect our country is on the first and leading place in the world. (Schmidt, 1941, quoted in Dubovitskaya, 2009, p. 150)
\eq
Here, “mathematical culture” is understood as a form of the intellectual capital of the nation, the total of mathematical knowledge, understanding, skills---and, crucially, problem-solving ability, covering all mathematics, including the highest levels of mathematical research. Trickier is the meaning of “mathematical culture” when these words are applied to individual people, as in Yaglom’s Foreword to Choquet (1970), the Russian translation of Choquet (1964):
\bq The book by Choquet assumes that the reader has certain mathematical culture” (p. 8),
\eq
which means general awareness of mathematics beyond the standard school or university courses, and some general content-independent mathematical traits. This concept is of course age dependent and has a different meaning when applied to a secondary school rather than to a university student or to a research mathematician---but awareness of abstraction being used in mathematics and preparedness at least to try to handle abstract concepts is assumed at a relatively early age.

Esakov (1994) gives a tiny, but exceptionally important piece of evidence of Stalin’s attitude to mathematics. The text of the speech by Lysenko at the infamous session of the All-Union Academy of Agricultural Sciences in 1948 (where Soviet genetics was totally destroyed, VASKhNIL, 1948) was submitted to Stalin for approval. Stalin underlined Lysenko’s statement “Every science is rooted in class [by its nature]” and wrote in the margin: “Hah-hah-hah… And mathematics? And Darwinism?” So, by 1948 (and perhaps even earlier), Stalin did not believe in the class nature of mathematics. This was a victory for mathematicians and had a profound effect on the fate of mathematics in the Soviet Union. Not every area of science was so lucky.

\section*{Mathematicians in the Reform}

One of the first published proposals of a reform in mathematics education was Boltyansky, Vilenkin, and Yaglom (1959); it did not introduce set-theoretic concepts, but development of geometry on the basis of geometric transformations featured in it prominently.

The new reformist program of mathematics teaching in 4th to 10th grades (ages 10 to 17) was approved by the Ministry of Education in 1968 and published in Matematika v Shlole [Mathematics in School], the mass circulation journal for mathematics teachers (Program, 1968) with implementation starting from 1970 (not simultaneously in all grades).

A draft of the program was published earlier (Kolmogorov, Markushevich, \& Yaglom 1967); it was indicated in the text that Kolmogorov and Markushevich were working on arithmetic, algebra, and elements of analysis, and Yaglom on geometry. The draft program was developed with participation of Vladimir Boltyansky, Yuri Makarychev, Galina Maslova, Konstantin Neshkov, Alexey Semushin, Antonin Fetisov, and Aleksandr Shershevsky.

I have already highlighted the names of Vladimir Boltyansky, Andrey Kolmogorov, Aleksey Markushevich, Naum Vilenkin, Isaak Yaglom as the key players in the reform. They all were prominent and internationally renowned research mathematicians, perhaps with the exception of Isaak Yaglom, a highly competent mathematician who devoted more time to university teaching and to his (truly remarkable!) work on popularization of mathematics than to writing research papers. He was definitely a mathematician, not a mathematical educationalist; for mathematicians he was “one of us.” Kolmogorov was the undisputed leader of the group.

\section*{Andrey Kolmogorov}

It could be conjectured that the political battles of the 1930s shaped Kolmogorov’s political stance. Perhaps he felt personally responsible for selection, nurturing, educating, developing professional mathematicians for the State and the Nation. He repeatedly called for increase in training of professional research mathematicians:
\bq
The Soviet Union nowadays needs large number of independent researchers of theoretical questions of mathematics. (Kolmogorov, 1959, p. 6)

Our country needs a large number of well–trained and talented mathematicians. It is very important that the professional mathematicians are chosen from those representatives of our youth who can work most productively in this area. One way of attracting gifted youth to mathematics is mathematical olympiads. Participation in school math circles and olympiads may help everyone to evaluate their own ability, seriousness, and strength of their passion for mathematics. (From the Editor, in Vasiliev and Egorov, 1963, p. 1)
\eq
These two texts were reprinted in Kolmogorov (1988).

Even more interesting is a letter from Kolmogorov to the psychologist Vadim Krutetskii (Kolmogorov, 2001) with comments on (the Russian original, published in 1968, of) Krutetskii’s famous book \emph{The Psychology of Mathematical Abilities in Schoolchildren} (Krutetskii, 1976). Kolmogorov gave Krutetskii advice on the use of statistics and also suggested how useful it could be to develop psychometric tests which allowed an early detection of mathematical abilities in children and also predicted the ceiling of their further development as mathematicians. It was obvious that Kolmogorov was interested only in one end product of mathematics education: professional research mathematicians. In the letter, Kolmogorov also remarked (I think on the basis of his work with kids in the Olympiad Stream, in particular, teaching in the specialist mathematics boarding school which he founded):
\bq
For now, as a practitioner, I am inclined to think that the nature of mathematical development achieved in accordance with the most modern recipes of early studies in set theory and algebra, up to the age of 10–13 years, and with fairly good success, could be replaced by a general education of quick wit and mental activity. But a delay in mastering strict logic and special mathematical skills at the age of 14–15 is becoming difficult to compensate.
\eq
But in his textbooks for the reform, he pushed the concept of equivalence relation and equivalence classes on 7th grade students (that is, aged 13) in mainstream education---we shall return to that later in this chapter.

Kolmogorov invested colossal amount of his time and effort in his reform. The bibliography in Shiryaev (2003) lists unbelievably 58 editions of school textbooks which he co-authored for the reform (and quite often, editions of the same title differed significantly from the previous ones) and 55 papers on various matters of the reform in Mathematics in School, a mass circulation journal for school teachers. This was really his reform, he owned it.

\section*{The Olympiad Stream}

It is not the aim of this chapter to give a detailed history of what we termed the Olympiad Stream, I will give only a summary description of the state it reached at 1970, the first year of the implementation of the reform, based in part on personal experience with occasional references to contemporary and modern sources.

Kolmogorov and his comrades-in-arms created a new subculture which valued advanced level elementary mathematics with focus on qualitative analytic thinking. They also created a community of people who shared these cultural values. “Mathematical culture” mentioned before was the culture of this community. The social aspects of this phenomenon are captured well by the title of the paper by Gerovich (2020): “Mathematical Paradise”: A Parallel Social Infrastructure of Postwar Soviet Mathematics.

By the 1970s, the Olympiad Stream had reached a considerable degree of development. First of all, it was a loose informal network of academic mathematicians and school teachers of mathematics involved in organizing mathematical competitions, running mathematical circles, summer schools, Sunday schools, distance learning by correspondence schools, etc.; of undergraduate students who helped to run all these activities; and, of course, of schoolchildren who were enthusiastically taking part in them. Next, it involved competitions and olympiads at every level: School, district/town, city/region, republican, national, as well as open national level olympiad by correspondence; in many cities, “mathematical battles” between schools. Boltyansky and Leman (1965) gives some idea of the remarkable mathematical quality of the top layer of these competitions. There was also a plethora of mathematical circles of various kinds and levels---at many schools, but also at universities and at cultural centres for children (the so-called “Houses of Young Pioneers”). Most circles were run by unpaid volunteers. In some cities, there were mathematical Saturday schools, Sunday schools, winter schools during the winter vacations in January, summer schools---the latter were usually run somewhere in countryside (Kolmogorov, Zhurbenko, Pukhova, Smirnova, \& Smirnov, 1971). Neretin (2021b) calls this plethora of activities the Konstantinov System, in memory of the great mathematical educator Nikolay Konstantinov (1932–2021), who was the principal contributor to its development.

This was a community with its own folklore, not always properly documented. Oral mathematical problems for various kinds of oral examinations, selection interviews, for use in “mathematical battles” featured in it prominently.

Moving to a more regular part of the systems, run by secondary schools, we have to mention non-compulsory facultative courses on additional chapters of mathematics at schools---they were supported by special textbooks, for example (Skopets, 1971), or a chapter in a textbook was reserved for facultative studies. An interesting example of the latter from the times of Kolmogorov’s reform was a chapter Logical Structure of Geometry (which included axiomatics for the planar geometry) in the geometry textbook for grade 8 (Kolmogorov, Semenovich, Gusev, \& Cherkasov, 1976).

Every significant town had specialist mathematics classes in some of the schools (a famous example is School 57 in Moscow (Sergeev, 2008), many cities had specialist mathematics schools.

A significant and indispensable role was played by mathematics and physics schools by correspondence: In the huge country, they reached everywhere. The three most well-known were run by the Moscow University---his was the first one, founded by the great mathematician Israel Gelfand in the early 1960s, and it set benchmarks for other schools (Rozov, Glagoleva, \& Rabbot, 1973). Another was run by the Moscow Physical-Technical Institute (Novoselov, 2010; Yumashev, 2012), and the last but not the least---the correspondence school of the Novosibirsk University (Khukhro, 2013). At one point, when I was a boy in Siberia, I was enrolled in all three---and this means, first of all, that information about them somehow reached me. No-one in my home school knew about that---actually many substreams of the Olympiad Stream were completely separate from the official school system.

Quite remarkable was the policy of publishing translations of books on mathematics, with a careful choice of the best available in the world literature. For example, easily available were Russian translations of such books as Choquet (1964), Courant and Robbins (1941), Coxeter (1961), Dieudonné (1964), Faure, Kaufmann, and Denis-Papin (1964), Hartshorne (1967), Niven (1961), Pólya (1962), Rademacher and Toeplitz (1957). In the Soviet academic publishing, there was a specific role: Translation editor, who supervised translation of a particular book and usually wrote a foreword. For Courant and Robbins and for Faure et al., the translation editor was Kolmogorov, for the six other books that I listed---Yaglom. Forewords by Kolmogorov and Yaglom put these books in the context of the new role of mathematics in science and in society---and hinted at the need to develop mathematical culture and mathematics education. I doubt that these books were popular among the majority of school teachers---they were published for people within the Olympiad Stream.

The Olympiad Stream itself produced a steady flow of excellent little books for children, often of high mathematical and didactic quality.\footnote{The full list of these books can be found on the site \href{https://math.ru/lib/}{https://math.ru/lib/} run by the
Moscow Centre for Continuous Mathematical Education, and on \href{https://mccme.ru/}{https://mccme.ru/},
\href{https://mccme.ru/index-e1.html}{https://mccme.ru/index-e1.html} which continues the proud traditions of the Olympiad Stream.} I wish to mention here four booklets, which were produced as assignments for the correspondence school at the Moscow University: Gelfand, Glagoleva, and Kirillov (1968), Gelfand, Glagoleva, and Shnol (1968), Kirillov (1970), Vasiliev and Gutenmakher (1970). They were examples of explaining mathematics to children in the simplest possible way---but without losing the essence of mathematics. For the mathematically experienced adult reader they were masterclasses of didactic transformation---I will say a few words about that in the next section.

Specialist mathematics schools require a special mention. More information on this particular phenomenon can be found in the paper Gerovich (2019) with the precisely chosen title “\emph{We Teach Them to Be Free}.” In addition to existing high-quality day schools in big cities (first of all, Moscow and Leningrad), four specialist boarding schools were set by a special decree of the Council of Ministers of the USSR, to be run by Moscow (Kolmogorov, Vavilov, \& Tropin 1981), Leningrad, Novosibirsk (Borovik, 2012), and Kiev Universities. The boarding school in Moscow became known as the Kolmogorov School, the one in Novosibirsk---as the Lavrentiev School, named after Mikhail Alekseevich Lavrentiev, an outstanding scholar who started his mathematical carrier as a student of Luzin (like Kolmogorov), but then turned to fluid dynamics and industrial mathematics (often with military applications). In late 1950s and 1960s, Lavrentiev founded the Siberian Branch of the Academy of Sciences, built an academic campus in the forest near Novosibirsk, and founded there the Novosibirsk University and the Physics and Mathematics Boarding School (PhMSh), my alma mater. I talk about my school in such detail because this sheds light on the main secret of the Olympiad Stream---it flourished because of the warm support from the Soviet military-industrial complex. Specialist mathematics schools continue to flourish in modern Russia (Konstantinov \& Semenov, 2021):
\bq
[S]chools with deeper study of mathematics (math schools) became the most important and very productive phenomenon in Russia’s education of the last decades. (p. 414)
\eq
And the last but not least---the Kvant magazine.``Kvant” means “Quantum” in Russian; it was a mass circulation monthly magazine on physics and mathematics for schoolchildren of grades 7 to 10 (but also devoted special pages to younger children). Kolmogorov was Kvant’s co-founder (in 1970) and the chief editor of mathematics.

This was a cultural stream which was recognized as such by most people who were actively involved in supporting it. The following are quotes from Voitishek (1973), the lecture notes of Vaclav Voitishek (1933–2003) given at the preparatory department of the Novosibirsk University:
\bq
The author assumes that the reader has access [...] to wonderful books about mathematical creativity [and refers to Russian translations of Pólya (1962) and Rademacher and Toeplitz (1957)]. (p. 3)
\eq
\bq
For the author, the models of exemplary exposition of mathematical truths are the books by Markushevich, Sikorsky, and Cherkasov (1967) on algebra, by Pogorelov (1972) on geometry, and by Boltyansky, Sidorov, and Shabunin (1972) on mathematics. (p. 6)
\eq

\section*{Didactic transformation}

The theoretical concept of didactic transformation from the mathematics education theory could be useful for explaining the principal reason of the failure of the reform. A compact formulation of what makes mathematics education so special can be found in a paper by Hyman Bass (2005):
\bq
Upon his retirement in 1990 as president of the International Commission on Mathematical Instruction, Jean-Pierre Kahane described the connection between mathematics and mathematics education in the following terms:
\bi
    \item In no other living science is the part of presentation, of the transformation of disciplinary knowledge to knowledge as it is to be taught (transformation didactique) so important at a research level.
    \item In no other discipline, however, is the distance between the taught and the new so large.
    \item In no other science has teaching and learning such social importance.
    \item In no other science is there such an old tradition of scientists’ commitment to educational questions. (p. 417)
    \ei
\eq
The concept of didactic transformation is fairly old and can be traced back to Auguste Comte (1852):
\bq
A discourse, then, which is in the full sense didactic, ought to differ essentially from one simply logical, in which the thinker freely follows his own course, paying no attention to the natural conditions of all communication. […] On the other hand, this transformation for the purposes of teaching is only practicable where the doctrines are sufficiently worked out for us to be able to distinctly compare the different methods of expanding them as a whole and to easily foresee the objections which they will naturally elicit. (Preface)
\eq
This concept is virtually unknown in English or Russian mathematics education literature (but apparently well known in France). It should be noticed that simple conversion of content in a “psychologically acceptable form” sometimes is not enough---a more serious mathematical work may be needed, and I will show you an example in the next section, where we shall try to apply this concept to assessment of reformist textbooks.

\section*{A First Case Study: Vectors}

In the ideal world, vectors (displacement, velocity, acceleration, force, etc.) are best introduced in physics courses. This was the approach in pre-reform Soviet schools, where vectors were only briefly mentioned in mathematics classes. This was the way of teaching in my years at the specialist boarding school at Novosibirsk by our physics lecturer Evgeny Bichenkov (lecture notes were later turned into a cute little book Bichenkov, 1999), and his approach was borrowed from the famous Feynman’s Lectures in Physics (Feynman, Leyton, \& Sands 1964). This was why my fellow-students and I were puzzled by Kolmogorov’s definition of vectors as parallel translations.
Boltyansky and Yaglom (1963), a few years prior to the reform, discussed equivalent forms of vector definitions:
Of course, no matter what definition we take, a vector from the elementary geometry’s point of view is a geometric object characterized by direction […] and length. But this definition is excessively general and does not trigger any geometric images. According to this general definition, a parallel translation is a vector […] since a parallel translation is characterized by its direction and length. Indeed, we could accept a definition: “vectors is any parallel translation.” This definition is logically perfect, and could be taken as a basis for development of the entire theory of actions over vectors and their applications. However this definition, despite its full correctness, also cannot satisfy us since thinking about vector as a geometric transformation appears to be insufficiently intuitive, distant from physical interpretations of vector magnitudes. (pp. 293–294, italics in the original)

This was an interesting warning; to show that it was justified we reproduce here a definition of parallel translation given by Lobeeva (1963) in a contemporary discussion of vectors in school education run by the magazine Mathematics in School:
\bq
A parallel translation is a point transformation of the plane, in which points trace equal, parallel, and equally directed segments. (p. 64)
\eq
It is a psychologically convincing description of a parallel translation as a process developing in time: moving points leave behind their traces showing their positions at intermediate moments of time. The trouble starts as soon as we consider other geometric transformations, for example rotations (which also can be seen as processes developing in time, so points are leaving traces behind them) and axial symmetries (where points just jump instantly across the axis---what are their traces?). Composition of two rotations through the same angle of 180 degrees, but around two different points is a parallel translation (this is easy to prove). But look at the trace of a point: It is the union of two half-circles. Composition of two axial symmetries in two parallel axes is also a parallel translation (this is even easier)---but where do traces of points come from?\footnote{An exercise for the reader: what is the composition of two axial symmetries with intersecting axes?}

This is an example of a didactic transformation gone wrong---because a concept was handled ignoring the wider mathematical context. In more specific terms, an element of the group of isometries of the Euclidean plane was treated by Lobeeva on its own, ignoring the rest of the group. Of course, Boltyansky and Yaglom knew this group and understood difficulties arising from defining vectors as parallel translations.
Verner (2012) explains that in the reform, the cautious approach of Boltyansky and Yaglom was overruled by Kolmogorov:
\bq
A. N. Kolmogorov volunteered to write a [geometry] textbook for the grades 6–8 forms. [\dots He] did not entrust writing “Geometry 6–8” to the well-known geometers V. G. Boltyansky and I. M. Yaglom who were in his commission. (pp. 19–20)
\eq
Kolmogorov preferred to stick to the definition of vector as a parallel translation. The wording in Kolmogorov, Nagibin, Semenovich, and Cherkasov (1977) was very casual, almost off the cuff:
\bq
In this chapter, we shall specifically deal with parallel translations, calling them by a new name: vectors. (p. 59, bold in the original)
\eq
Unfortunately, the definition of parallel translation starts in Kolmogorov, Semenovich, and Cherkasov (1979) with a definition of equivalence relation and then immediately states, without proof, that an equivalence relation on a set partitions this set into equivalence classes. This is followed by this example:
\bq
The parallelity relation between straight lines in the plane defines a partition of the set of straight lines in the plane into classes. Each of these classes consists of straight lines, parallel to each other […]. These classes are bundles of parallel lines. Another example of these classes is directions… (p. 128, italics in the original)
\eq
No other examples are given, and a definition of direction appears later:
\bq
The set of rays, each of which is co-directed with the same ray, is called direction. (p. 130, italics in the original)
\eq
Finally, a definition of parallel translation is given on p. 132---quite similar in wording to the one given in Klopsky, Skopets, and Yagodovsky (1980), mentioned in the Introduction and quoted by Pontryagin in his famous attack on the reform (Pontryagin, 1980).

It appears that at least some of Kolmogorov’s collaborators understood the difficulty of his approach. In a book for teachers by Gusev, Kolyagin, and Lukankin (1976), five pages (pp. 6–11) are devoted to explaining, to teachers, Kolmogorov’s definitions and its versions used in various textbooks. The equivalence relation features prominently (p. 8). An advanced set-theoretic approach is also invoked:
\bq
To summarize, we considered a possibility of introducing the concept of vector as a set of pairs of points defining the same parallel translation, that is, the set of all pairs (X, Y), for which T(X) = Y, as a vector. The set of pairs (X, Y) is sometimes called the graph of the parallel translation.

In the modern treatment it is conventional to identify the graph of the mapping with the mapping itself. Everything said before has led to the identification, in the school course of mathematics, of a parallel translation and a vector as synonyms, denoting exactly the same concept. (p. x)
\eq
In the reform, the set theory was confined to non-compulsory enhancement courses---but I was unable to find anywhere in them the delicate and abstract identification of a function (that is, a map from a set to a set) with its graph. A proper definition of the graph required a definition of the direct product of two sets---which was also nowhere to find.

	Nowadays in England, my university colleagues consider the theorem about an equivalence relation partitioning a set into equivalence classes as pons asinorum of undergraduate abstract mathematics. Alas, many graduates from English universities obtain their bachelor’s degree in mathematics without grasping this concept. At that time, circa 1970 in Russia, the equivalence/partition duality was perhaps one of the boundary markers between the Olympiad Stream and the mainstream school mathematics. Kolmogorov himself, and the reformer mathematicians of his circle were experts in education of kids in the Olympiad Stream. But they crossed the boundary into mainstream education without caring about didactic transformation of the new material which they brought with themselves.

In short, at the methodological level, the principal reason for failure of the reform was absence (or failure) of appropriate didactic transformation of the new mathematical content. This was even more surprising because inside of the Olympiad Stream, didactic transformation was used quite successfully---I have already given a few examples.

The famous geometer Aleksandr D. Aleksandrov (1980, reprinted as Aleksandrov, 2008a) gave a rather harsh assessment of the new course of geometry:
\bq
It is hard to find something more harmful for the spiritual---mental and moral---development, than train a person to pronounce words which meaning he does not really understand, and, when necessary, is guided by other concepts. (p. 307)
\eq
It is interesting that in his speech at the meeting of the Academic Council of the Institute of Mathematics of the Siberian Branch of the Academy of Sciences USSR on December 25, 1980 (printed as Alexandrov 2008b), he harshly criticized the paper Pontryagin (1980) as an attack on “abstract” mathematics and put the blame for the botched reform on Kolmogorov’s collaborators (Zalman Skopets was the only one named---perhaps because he was not a mathematician, but a mathematical educationalist).
\bq
We talk not about a set-theoretic approach, not about some special abstractions and sophistry, but about very simple things, like crude mistakes in Russian language in Geometry 6 or a ridiculous definition of a polytope in a textbook for grades 9–10. It is not abstraction in mathematics, but, in the final count, abstracting from responsibility, abstracting from conscientiousness are the root of mistakes and absurdities in school teaching as well as in public pronouncements about mathematics. (p. 319)
\eq

\section*{A Second Case Study: Probability Theory}

Kolmogorov was one of the founders of the modern probability theory. His ground-breaking work \emph{Grundbegriffe der Wahrscheinlichkeitsrechnung} [Basic concepts of probability theory] (Kolmogoroff, 1933) was a masterpiece of exposition of mathematics.

So it could be surprising that in the reform, probability theory (and some elementary combinatorics needed) were limited in scope and confined to a short chapter in (non-compulsory) facultative courses in grade 9. The theory was restricted to the finite frequentist setting and Bernoulli trials, it just barely touched conditional probability, and a lot of attention was given to the direct computation of probabilities with the help of combinatorial formulae; however, it included some simple examples of geometric probabilities (Kotii \& Potapov, 1971).

This modest original treatment was soon developed to include random variables and one of the simplest versions, due to Chebyshev, of the law of large numbers (Antipov, Vilenkin, Ivashev-Musatov, \& Mordkovich, 1979; Firsov, Bokovnev, \& Shvartsburd, 1977). Still the content appeared to be rather unimaginative. There was no sign of a contribution from Kolmogorov to the probability theory chapters of textbooks. Why? Perhaps, being the expert in probability theory, he understood that any step away from the elementary---and frequently artificial---material led into serious conceptual difficulties.

David Corfield, a mathematician and well-known philosopher of mathematics, made the following incisive comments in the context of debates around teaching probability and statistics in schools in England (D. Corfield, personal communication, October 12 and 13, 2010):
\bq
	One intriguing problem about teaching probability theory is that there are at least four distinct interpretations of probability (an objective and a subjective Bayesianism, a propensity theory, and a frequency theory), along with various pluralist positions. Unless you work in artificial situations with, say, perfect dice, these differences, which I imagine most school teachers are unaware of, will confuse one’s teaching.

Presumably an analysis of decisions to play lotteries could be done in a fairly uncontroversial way, though the relative utility of losing a pound and gaining so many millions is far from obvious.

Then there’s the question of various forms of optional insurance. Here we enter the problems of assessing likelihoods of events when data only covers certain groups. E.g., how do I calculate the probability of suffering a heart attack when all I have is data for, say, non-smoking 40-year-old males. Maybe there’s no data for those with my diet, exercise, income, job satisfaction, marital happiness, etc. Ultimately there’s only one me. This is the “reference class” problem.

Odds in horse races provide a very good illustration of probabilities. Are they
\bi
            \item[a)] the unique propensity of a horse in that precise situation to win;

            \item[b)] the limiting frequency in some long series of events;

            \item[c)] a measure of subjective expectation, reflected in betting behavior;

            \item[d)] an objective measure of the expectations of a rational agent given certain information?
\ei
            \eq
For a more detailed discussion, see Gillies (2000).

Perhaps we have to conclude that Kolmogorov treated probability theory very differently from geometry.

\section*{Social Blindness}

Insufficient attention to didactic transformation of new material, making it accessible if not to all, but at least to the majority of mainstream students may be linked to reformers’ surprising social blindness in many other aspects of the reform.

For example, the need to re-educate, retrain the whole army of teachers was somehow overlooked. I made a systematic search for books for teachers and students in pedagogical colleges which were produced in the support of the reform. Not much; it was not surprising that Kolyagin (2001) characterized the provision of advice to teachers, of didactic materials etc. as essentially non-existent. Kolyagin also pointed to the disruption of teaching the new generation of teachers in pedagogical colleges---something that could be foreseen before the launch of the reform and appropriately alleviated.

As an unexpected side effect, the neglect of teachers damaged the Olympiad Stream. In 1974–1978, I was a university students, but I was involved in running regional mathematics olympiads in Siberia and Soviet Far East and had a chance to see the negative effect of the reform on the participants: Their overwhelmed teachers could no longer give them enough attention and time.

We provide another example of blindness to social realities. In their experimental textbook Boltyansky, Volovich, and Semushin (1979) gave the following advice to 12-year-old students (grade 6):
\bq
You will not find answers to the problems at the end of the book. After all, we want you to learn to reason in the right way, to have confidence in the correctness of the logical arguments, in the rationality of the solution found. We want you to be able to apply your geometric and logical knowledge in life, in your future work---but life does not provide answers to questions, easy or difficult, that it poses. Therefore, get used to solving problems without “peeping” in solutions. And if sometimes your solution turns out to be not entirely correct, imprecise, you will be corrected by your comrades and the teacher. (pp. 3–4)
\eq
Alas, after five years of attempts to learn mathematics most student were well aware that no-one would correct them, if answers had not been given to their teachers in advance. Advising students to seek help from teachers was a breach of one of the unwritten traditional rules for mathematics textbooks in Russia: A good student should be able to learn mathematics directly from a textbook without help from a teacher. However, within the Olympiad Stream it was fine not to give answers to problems in mathematics and physics; moreover, it was a normal practice. I myself was taught in this tradition; for example, Bichenkov (1999), the book which grew up from Bichenkov’s lectures in Novosibirsk PhMSh which I attended--- had no answers. Problems were wonderful---but I still do know how to solve many of them.
And the last example has some curiosity value. For the promotion of the reform, Boltyansky and Levitas (1973) tried to appeal directly to parents and wrote a book in the form of a dialogue between a mathematician and a few parents. Alas, the parents in their book are not very representative of the general population. In the book, they (and the reader) are offered “homework.” Here is one of the problems:
On return from school, your son asked you a question:
\bq
We have been told that axioms cannot be proven, and gave an axiom: “There is only one straight line that passes through two points.” Why couldn’t this be proven? Apply the ruler and check that the second line goes along the first one. This is a proof!
What would you say to your son? (p. x)
\eq
How many real parents, even after reading the book, were able to coherently answer this question?

I have a feeling that the reformers were not aware of many aspects of the socio-economic situation in the country even those which directly affected education. Perhaps they sincerely believed in the official dogma of the social homogeneity of the Soviet society. This was an illusion. In the next section, I will try to explain the roots of this illusion.

\section*{The Golden Age of Soviet Mathematics Education}

This was the principal mistake of the reformers: They continued to live in the Golden Age of Soviet mathematics, as the period from the 1950s through the 1970s is frequently called (Gerovitch, 2013, 2019, 2020), and they did not realize that it was coming to an end. I was lucky to see the last days of the Golden Age, and I have seen how it ended. And I cannot blame them for their illusion.

In the Soviet Union, the society was not almost homogeneous, as the official propaganda insisted, it was deeply stratified, and had an etacratic social structure---the social position of people was almost entirely determined by their place within the dictatorial autocratic state and the centrally planned economy (Radaev, \& Shkaratan, 1995; Shkaratan, 2012).

In the beginning of the 1970s, all aspects of life in the Soviet Union were still dominated by the tidal wave of social mobility unprecedented in history of humanity. Social mobility ensured the stability of the society and of the totalitarian system that ruled the society. No matter how hard life was, every family could have realistic expectations that their children would have better lives if they got education---and, for that reason, people were prepared to forgive the blunders and even crimes of their rulers. But---and this an important but---the social mobility was something that was centrally planned, as everything in the economy. It was planned, for example, how many people will leave villages and move to cities to work in industry.

It was a plain numerical horizontal expansion of the economy, not accompanied by the growth in productivity. In the case of education, especially mathematics education, the expansion could be seen with exceptional clarity. Expanding school education required more school teachers, more pedagogical colleges and universities, more university teachers, with demand feeding back into the need for further improvement and expansion of schools. It was driven by an insatiable, it looked at the time, demand from the military-industrial complex for educated (and, above all, mathematically literate) workers and engineers, and from the armed forces for soldiers and officers with good mathematical skills.

In more politically sensitive areas, the process of social mobility was skillfully manipulated. The authorities were running an elaborate system of positive discrimination and promotion of young people from politically safe strata of the society---and this involved institutionalized anti-Semitism.

Of course, the economic expansion without a matching growth in productivity cannot last, and, exactly in the 1970s, the music stopped. For the social mobility, there was no more room at the top. The Soviet Union collapsed very soon afterwards.

It can be conjectured that the reformers were either blinded by the official propaganda, or forced to behave as if they believed it.
The limited space of this chapter does not allow me to provide an in-depth analysis of my observations---I hope this will eventually be done by professional historians and politologists. My role was simple: I wished to attract these experts’ attention to an exciting object of study.

\section*{A Lesson for our Times?}

The principal lesson of the Kolmogorov’s reform: The methods of the Olympiad Stream are not transferable to mainstream education on the cheap, as it was attempted by the reformers. What was needed for success was serious investment in a proper reform, with sociological studies, with more time spent on developing and testing textbooks, books for teachers, didactic material, textbooks for pedagogical colleges, with systematic re-education and professional development of teachers, with a much larger and better coordinated team of developers, and proper project management. In the economic situation of the Soviet Union of 1968 all that was unfeasible.

In conclusion, I offer to the reader’s attention a short fragment from a blatantly self-promotional film from Yandex, the Russian IT and Internet giant (Yandex 2020)---it has English subtitles. I recommend to watch just the segment 3:22–4:35, where Tigran Khudaverdyan, General Director and Director of Operations at Yandex, answers a question from a reporter: “What is yandexoid?” The word “yandexoid” entered Russian IT jargon, it is the proud self-designation of Yandex employees who apparently feel themselves being the salt of the earth (Matthew 5:13).
\bq
Alexei Pivovarov (the reporter): What is yandexoid?\\
Tigran Khudaverdyan: yandexoid\dots\ yandexoid is\dots\ it's such an environment\dots\ Lesha, have you ever participated in olympiads, some school ones?\\
Alexei Pivovarov: There was a case. It didn’t end well.\\
Tigran Khudaverdyan: My memories are: You are the first guy in the village, in your class, in your school, such a great fellow. You think you know the topic best. You come there, for example, to the city olympiad\dots\ but it doesn't matter where, from the city one to the republican one, and you understand that in general everyone there is smarter than you, that all your greatness is broken, just by the fact that the strongest have gathered there. Well, to be an yandexoid is to be an olympiadnik every day [italics added]. You come to work every day, you have to be able to be wrong, not knowing, and then you came up with something absolutely brilliant. They may show you very reasonably that it won't work there---or vice versa. And if you can't stand that, then you won't be able to work.
\eq
This is the new reality: The selective Olympiad Stream, with its intrinsic competitiveness, is welcome in the corporate world.

Meanwhile, nowadays in Russia, in the new technological and socio-economic environment, the debate about the balance between the selective stream and mainstream, and the content of the mainstream in school mathematics education is very much alive (Borovik, Kocsis, \& Kondratiev, 2022; Khalin, Vavilov, \& Yurkov, 2022; Konstantinov \& Semenov, 2021). It is important to ensure that it is informed by lessons from Kolmogorov’s reform.

\section*{Acknowledgements}

The author thanks Glen Aikehead, Szabó Csaba, Gregory Cherlin, Roman Kossak, Dmitrii Pasechnik, Nikolai Vavilov, Alexander Veselov, and Theodore Voronov for their feedback and advice. He thanks Dirk De Bock for his patient and detailed editing of this chapter.

\section*{References}

Abramov, A. M. (2016). \emph{The great world of Fatherland, or the Kolmogorov project for the XXI century} [in Russian]. Sankt-Peterburg, Russia: Obrazovatel'nye proekty

Aleksandrov, A. D. (1980). On geometry [in Russian]. \emph{Matematika v Shkole}, 3, 56--62.

Aleksandrov, A. D. (2008a). On geometry in school [in Russian]. In \emph{Selected works}. Vol. 3: \emph{Articles from various years} (pp. 296--308). Novosibirsk, Russia: Institut matematiki im. S. L. Soboleva.

Aleksandrov, A. D. (2008b). On the state of school mathematics [in Russian]. In \emph{Selected works}. Vol. 3: \emph{Articles from various years} (pp. 309--325). Novosibirsk, Russia: Institut matematiki im. S. L. Soboleva.

Antipov, I. N., Vilenkin, N. Ya., Ivashev-Musatow, O. S., \& Mordkovich, A. G. (1979). \emph{Selected questions of mathematics. Grade 9. A facultative course} [in Russian].

Bass, H. (2005) Mathematics, mathematicians, and mathematics education.\emph{ Bulletin of the American Mathematical Society}, 42(4), 417--430.

Bichenkov, E. I. (1999). \emph{Laws of mechanics} [in Russian]. Novosibirsk, Russia: Izdatel’stvo IDMI.

Boltyansky, V. G., \& Leman, A. A. (1965). \emph{Collection of problems of Moscow mathematical olympiads} [in Russian]. Moscow, Russia: Prosveshchenie.

Boltyansky, V. G., \& Levitas, G. G. (1973). \emph{Mathematics attacks parents} [in Russian]. Moscow, Russia: Pedagogika.

Boltyansky, V. G., Sidorov, Yu. V., \& Shabunin, M. I. (1972). \emph{Lectures and problems on elementary mathematics }[in Russian]. Moscow, Russia: Nauka.

Boltyansky, V. G., Vilenkin, N. Ya., \& Yaglom, I. M. (1959). About the content of the course of mathematics in secondary school [in Russian]. \emph{Mathematics, its Teaching, Applications, and History}, 4, 131--143.

Boltyansky, V. G., Volovich, M. B., \& Semushin, A. D. (1979). \emph{Geometry. A trial textbook for grades 6--8} [in Russian]. Moscow, Russia: Prosveshchenie.

Boltyansky, V. G., \& Yaglom, I. M. (1963). Vectors and their uses in geometry [in Russian]. In P. S. Alexandrov \& A. I. Markushevich (Eds.), \emph{Encyclopedia of elementary mathematics. Book IV: Geometry} (pp. 292–381). Moscow, Russia: GIFML.

Boltyansky, V. G., \& Yaglom, I. M. (1965). The school mathematical circle at the Moscow State University and Moscow mathematical olympiads [in Russian]. In V. G. Boltyansky \& A. A. Leman (Eds.),\emph{ Collection of problems of Moscow mathematical Olympiads }(pp. 3--46). Moscow, Russia: Prosveshchenie.

Bonchkovsky, R. N. (1936). \emph{Moscow mathematical olympiads of 1935 and 1936 years} [in Russian]. Moscow, Russia: ONTI NKTP SSSR.

Borovik, A. V. (2012). “Free Maths Schools”: Some international parallels. \emph{The De Morgan Journal} 2(2), 23--35. Retrieved May 23, 2021, from\linebreak \href{https://tinyurl.com/355ac33c}{https://tinyurl.com/355ac33c}.

Borovik, A. V., Karakozov, S. D., \& Polikarpov, S. A. (2021). \emph{Mathematics education policy as a high stakes political struggle: The case of Soviet Russia of the 1930s.} arXiv:2105.10979 [math. HO]. Retrieved May 23, 2021, from \href{https://arxiv.org/abs/2105.10979}{https://arxiv.org/abs/2105.10979}.

Borovik, A., Kocsis, Z., \& Kondratiev, V. (2022). Mathematics and mathematics education in the 21st century.  [To appear.]

Chebotarev, N. G. (1937). Mathematics olympiad of schoolchildren in Kazan [in Russian]. \emph{Matematicheskoe Prosveshchenie}, 11, 65.

Chistyakov, I. I. (1935). The mathematical olympiad of the Leningrad State University named after A. S. Bubnov [in Russian]. \emph{Matematiheskoe Prosveshchenie}, 3, 59--63.

Choquet, G. (1964). \emph{L’enseignement de la géométrie} [Teaching of geometry]. Paris, France: Hermann.

Choquet, G. (1970). \emph{Geometry} [in Russian]. Moscow, Russia: Mir.

Comte, A. (1852). \emph{Catéchisme positiviste, ou, Sommaire exposition de la religion universelle, en onze entretiens systématiques entre une femme et un Prêtre de l’humanité} [Positivist catechism, or, Summary exposition of universal religion, in eleven systematic interviews between a woman and a Priest of humanity]. Paris, France: Author.

Courant, R., \& Robbins, H. (1941). \emph{What is mathematics?} London, United Kigdom---New York, NY---Toronto, Canada: Oxford University Press.

Coxeter, H. S. M. (1961). \emph{Introduction to geometry}. New York, NY---London, United Kingdom: John Wiley \& Sons, Inc.

Demidov, S. S., \& Levshin, B. V. (Eds.). (1999). \emph{The case of academician Nikolai Nikolaevich Luzin} [in Russian]. Moscow, Russia: RkhGI.

Dieudonné, J. (1964). \emph{Algèbre linéaire et géométrie élémentaire} [Linear algebra and elementary geometry]. Paris, France: Hermann.

Dubovitskaya, M. A. (2009). The activities of O.Yu. Schmidt in the Moscow University [in Russian]. \emph{Historico-Mathematical Investigations}, 13(48), 138--153.

Esakov, V. D. (1994). New on the session of VASKhNIL of 1948 [in Russian]. In M. G. Yaroshevsky (Ed.), \emph{Repressed science} (Issue II, pp. 57--75). St. Petersburg, Russia: Nauka.

Faure, R., Kaufmann, A., \& Denis-Papin, M. (1964). \emph{Mathématiques nouvelles.} Tome 1 [New mathematics. Vol. I]. Paris, France: Dunod.

Feynman, R. P., Leyton, R. B., \& Sands, M. (1964). \emph{Lectures on physics}. Vol. 1. Reading, MA: Addison Wesley.

Firsov, V. V., Bokovnev, O. A., \& Shvartsburd, S. I. (1977). \emph{The state and perspectives of facultative studies in mathematics. An aid for teachers} [in Russian]. Moscow, Russia: Prosveshchenie.

Fomin, D. (2020). Mathematical ``archeology'': Problems of the first Soviet school olympiad in mathematics [in Russian]. \emph{Kvant}, 7, 16--21.

Gelfand, I. M., Glagoleva, E. G., \& Kirillov, A. A. (1968). \emph{Method of coordinates} [in Russian]. Moscow, Russia: Nauka.

Gelfand, I. M., Glagoleva, E. G., \& Shnol, E. E. (1968). \emph{Functions and graphs} [in Russian]. Moscow, Russia: Nauka.

Gerovitch, S. (2013). Parallel worlds: Formal structures and informal mechanisms of postwar Soviet mathematics. \emph{Historia Scientiarum}, 3(22), 181–200.

Gerovitch, S. (2019). ``We teach them to be free''. Specialized math schools and the cultivation of the Soviet technical intelligentsia.
Kritika: \emph{Explorations in Russian and Eurasian History}, 4(20), 717--754.

Gerovitch, S. (2020). ``Mathematical paradise'': A parallel social infrastructure of postwar Soviet mathematics. \emph{Logos}, 2(30), 93--128.

Gillies, D. (2000). \emph{Philosophical theories of probability}. London, United Kingdom: Routledge.

Grigoriev, E. I. (1937). Problems for schoolchildren [in Russian]. \emph{Matematicheskoe Prosveshchenie}, 11, 67--68.

Gusev, V. A., Kolyagin Yu. M., \& Lukankin G. M. (1976). \emph{Vectors in the school course of geometry} [in Russian]. Moscow, Russia: Prosveshchenie.

Hartshorne, R. (1967). \emph{Foundations of projective geometry}. New York, NY: W. A. Benjamin, Inc.

Karp, A., \& Vogeli, B. (Eds.). (2010). \emph{Russian mathematics education: History and world significance}. Singapore: World Scientific.

Khalin, V., Vavilov, N., \& Yurkov, A. (2022). The skies are falling: Mathematics for non-mathematicians. \emph{Mathematics}. [To appear.]

Khukhro, E. I. (2013). Physics and mathematics school by correspondence at the Novosibirsk State University. \emph{The De Morgan Journal}, 1(3), 1--6. Retrieved May 23, 2021, from \href{https://tinyurl.com/bwnf2z4c}{https://tinyurl.com/bwnf2z4c}.

Kirillov, A. A. (1970). \emph{Limits} [in Russian]. Moscow, Russia: Nauka.

Klopsky, V. M., Skopets Z. A., \& Yagodovsky M. I. (1980). \emph{Geometry. A textbook for grades 9 and 10 of secondary school} [in Russian]. Moscow, Russia: Prosveshchenie.

Kolmogoroff, A. (1933). \emph{Grundbegriffe der Wahrscheinlichkeitsrechnung} [Basic concepts of probability theory]. Berlin, Germany: Springer.

Kolmogorov, A. N. (1959). \emph{On the profession of a mathematician} [in Russian]. Moscow, Russia: Izdatel'stvo Moskovskogo Universiteta.

Kolmogorov, A. N. (1988). \emph{Mathematics---science and profession} [in Russian]. Moscow, Russia: Nauka.

Kolmogorov, A. N. (2001). On development of mathematical abilities [Letter to V. A. Krutetskii, in Russian]. \emph{Voprosy Filosofii}, 3, 103--106.

Kolmogorov, A. N., Semenovich, A. F., \& Cherkasov, R. S. (1979). \emph{Geometry. Textbook for grades 6--8 of secondary school} [in Russian]. Moscow, Russia: Prosveshchenie.

Kolmogorov, A. N., Semenovich, A. F., Gusev, V. A., \& Cherkasov, R. S. (1976). \emph{Geometry. Textbook for grade 8 of secondary school} [in Russian]. Moscow, Russia: Prosveshchenie.

Kolmogorov, A. N., Semenovich, A. F., Nagibin, F. F., \& Cherkasov, R. S. (1977). \emph{Geometry. Textbook for grade 7 of secondary school} [in Russian]. Moscow, Russia: Prosveshchenie.

Kolmogorov, A. N., Vavilov, V. V., \& Tropin, I. T. (1981). \emph{Physics and mathematics school at MGU} [in Russian]. Moscow, Russia: Znanie.

Kolmogorov, A. N., Zhurbenko, I. G., Pukhova, G. V., Smirnova, O. S., \& Smirnov, S. V. (1971). \emph{Summer school at Rubskoi lake} [in Russian]. Moscow, Russia: Prosveshchenie.

Kolyagin, Yu. M. (1977). \emph{Problems in mathematics education. Part I: Mathematical problems as means of teaching and development of students. Part II: Teaching mathematics via problems and teaching problem solving} [in Russian]. Moscow, Russia: Prosveshchenie.

Kolyagin, Yu. M. (2001). \emph{Russian school and mathematics education. Our pride and our pain} [in Russian]. Moscow, Russia: Prosveshchenie.

Kolyagin, Yu. M., \& Savvina, O. A. (2012). \emph{Rebellion of the Russian ministry and Mathematics Division of Academy of Sciences USSR} [in Russian]. Elets, Russia: Elets State University.

Konstantinov, N. N., \& Semenov, A. L. (2021). Resultative education in mathematics schools [in Russian]. \emph{Chebyshevskii Sbornik}, 22(1), 413--436.

Kostenko, I. P. (2013). \emph{The problem of quality of mathematics education in the light of the historic perspective} [in Russian]. Moscow, Russia: ROSZhELDOR.

Kostenko, I. P. (2014). 1965--1970. Organisational preparation of reform-70: Ministry of Education, Academy of Pedagogical Sciences, appointments, programs, textbooks (article five) [in Russian]. \emph{Matematicheskoe Prosveshchenie}, 3(71), 2--18.

Kotii, O. A., \& Potapov, V. G. (1971). \emph{Elements of probability theory with elements of combinatorics} [in Russian]. In Z. A. Skopets (Ed.), Collection of problems on mathematics (pp. 97--116). Moscow, Russia: Prosveshchenie.

Krutetskii, V. A. (1976). \emph{The psychology of mathematical abilities in schoolchildren}. Chicago, IL: University of Chicago Press.

Lobeeva, A. A. (1963). On vectors in the school mathematics course [in Russian]. \emph{Mathematics in School}, 2, 64--68.

Markushevich, A. I., Sikorsky, K. P., \& Cherkasov, R. S. (1967). \emph{Algebra and elementary functions} [in Russian]. Moscow, Russia: Prosveshchenie.
	
Neretin, Yu. A. (2019). The Kolmogorov reform of mathematical education, 1970--1980 [in Russian]. Retrieved July 4, 2021, from\\ \href{https://arxiv.org/pdf/1911.06108.pdf}{https://arxiv.org/pdf/1911.06108.pdf}.

Neretin, Yu. A. (2021a). \emph{Time of Luzin. The birth of the Moscow mathematical school: Soviet mathematics at the background o social cataclysms of the 1920--1930s years} [in Russian]. Moscow, Russia: URSS.

Neretin, Yu. A. (2021b). Nikolay Konstantinov and the Konstantinov System [in Russian]. Retrieved November 7, 2021, from \href{https://arxiv.org/pdf/2110.03621.pdf}{https://arxiv.org/pdf/2110.03621.pdf}.

Niven, I. (1961). Numbers: Rational and irrational. New York, NY: Random House.

Novoselov, K. (2010). Konstantin Novoselov---Biographical. Nobel Prize Outreach AB. Retrieved September 18, 2021, from\\ \href{https://www.nobelprize.org/prizes/physics/2010/novoselov/biographical}{https://www.nobelprize.org/prizes/physics/2010/novoselov/biographical}.

Phillips, C. J. (2015). \emph{The New Math. A political history}. Chicago, IL: University of Chicago Press.

Pogorelov, A. V. (1972). \emph{Elementary geometry} [in Russian]. Moscow, Russia: Nauka.

Pólya, G. (1962). \emph{Mathematical discovery: on understanding, learning and teaching problem solving}. New York, NY: Wiley.

Pontryagin, L. (1980). About mathematics and quality of its teaching [in Russian]. \emph{Kommunist}, 14, 99--112.

Program. (1968). The program of mathematics for secondary school [in Russian]. \emph{Matematika v Shkole}, 2, 5--20.

Radaev, V. V., \& Shkaratan, O. I. (1995). \emph{Social stratification }[in Russian]. Moscow, Russia: Nauka.

Rademacher, H., \& Toeplitz, O. (1957). \emph{The enjoyment of mathematics. Selections from mathematics for the amateur}. Princeton, NJ: Princeton University Press.

Rozov, N. Kh., Glagoleva E. G., \& Rabbot, Zh. M. (\emph{1973). Mathematical correspondence school at MGU} [in Russian]. Moscow, Russia: Znanie.

Schmidt, O. Yu. (1941). On connection between science and practice [in Russian]. Archive RAS, fund 496, inventory 1, file 248.

Sergeev, P. W. (2008). \emph{Mathematics in specialist classes of School 57. Mathematical analysis} [in Russian]. Moscow, Russia: MtsNMO.

Sharygin, I. (2002). Which “horse” will bring death to Russian mathematics? [in Russian]. \emph{Otechestvennye Zapiski}, 3(2). Retrieved July 4, 2021, from \href{https://strana-oz.ru/2002/2/ot-kakogo-konya-primet-smert-rossiyskaya-matematika}{https://strana-oz.ru/2002/2/ot-kakogo-konya-primet-smert-rossiyskaya-matematika}.

Shiryaev, A. N. (Ed.). (2003). \emph{Kolmogorov. Jubilee publication in three books}. Book 1. \emph{Truth is blessing. Bibliography} [in Russian]. Moscow, Russia: Fizmatlit.

Shkaratan, O. I. (2012). Sociology of inequality. \emph{Theory and reality} [in Russian]. Moscow, Russia: Izdatelskii dom Vyshei shkoly ekonomiki.

Skopets, Z. A. (Ed.). (1971). \emph{Collection of problems on mathematics} [in Russian]. Moscow, Russia: Prosveshchenie.

Strunnikov, V. A., \& Shamin, A. N. (1989).\emph{ T. D. Lysenko and lysenkoism. Destruction of Soviet genetics in the 30–40s} [in Russian]. \emph{Biologia v Shkole}, 2, 15-–20.

Vasiliev, N. B., \& Egorov, A. A. (1963). \emph{Collection of preparatory problems for All-Russia olympiad of young mathematicians} [in Russian]. Moscow, Russia: Uchpedgiz.

Vasiliev, N. B., \& Gutenmakher, V. L. (1970). \emph{Straight lines and curves} [in Russian]. Moscow, Russia: Nauka.

VASKhNIL. (1948). \emph{On the situation in the biological science. Stenographic report of the session of All-Union Academy of Agricultural Sciences named after V.I. Lenin. July 31 to August 7, 1948} [in Russian]. Moscow, Russia: OGIZ-Sel’khozgiz.

Vavilov, N. A. (2022). Computers as the new reality of mathematics: IV. Goldbach’s problem. [In Russian]. \emph{Kompyuternye Instrumenty v Obrazovanii}. [To appear.]

Verner, A. L. (2012). A. D. Aleksandrov and the course of school mathematics [in Russian]. \emph{Matematicheskie Struktury i Modelirovanie}, 25, 18-–38.

VKP(b). (1934). \emph{XVII congress of VKP(b)}. [Stenographic report, in Russian]. Moscow, Russia: Politizdat. Retrieved July 4, 2021, from\\ \href{http://istmat.info/node/52145}{http://istmat.info/node/52145}.

Voitishek, V. V. (1973). \emph{Lectures on mathematics for students of the preparatory department} [in Russian]. Novosibirsk, Russia: Novosibirsk State University.

VSM. (1935). Resolutions of the Second All-Union Mathematical Congress [in Russian]. \emph{Matematicheskoe Prosveshchenie}, 3, 52--59.

Yandex (2020). YaC 2020: How we make Yandex. YouTube\\ \href{https://www.youtube.com/watch?v=qlFTxADGixA}{https://www.youtube.com/watch?v=qlFTxADGixA}, premiered November 25, 2020. Retrieved May 15, 2021.

Yumashev, D. (2012). ZFTSh: A specialist correspondence school. \emph{The De Morgan Journal}, 	2(2), 37--41. Retrieved May 23, 2021, from htps://tinyurl.com/7sc5zs6e.

Zalgaller, V. A. (2021). I was lucky to live surrounded by people who loved science [in 	Russian]. \emph{Matematicheskoe Prosveshchenie}, s. 3, 28, 15--34.

\end{document}